\RequirePackage{fix-cm}
\documentclass[12pt]{amsart}           

\usepackage{graphicx}
\usepackage{url}
\usepackage{hyperref}
\hypersetup{colorlinks=true,linkcolor=blue,urlcolor=blue}
\usepackage{mathptmx}      
\usepackage{amssymb}
\newcommand{\bq}{\begin{quote}}
\newcommand{\eq}{\end{quote}}

\begin{document}
\begin{abstract}
University level mathematics in a number of countries is under pressure to `decolonise
the curriculum'. This paper considers, as a test case, a possible `decolonisation' of
linear algebra.  This is a representative case, since linear algebra is one of the core
courses of undergraduate mathematics and a mathematical discipline with a millennia
long historic tradition. This paper is written for my colleagues, university
mathematicians. In my opinion, it could help them to determine their position and calmly  stick
to it without entering into an unnecessary debate with promoters of `decolonisation of
the curricula'. The paper offers  a simple and honest defence against
`decolonisation' pressures: tell students the real (pre)history of a particular
mathematical discipline. Let us call this activity `historical enrichment'.

It would be useful if all attempts at `decolonization' (and vice versa, `historical
enrichment') were known to a wider circle of the mathematical community. Publicity
and an open discussion are the best way to resist outside pressures to engage in virtue
signalling at the expense of historical and mathematical truth.

The international mathematical community should defend academic freedom and insist
on our right to formulate our curricula and evaluate the history of mathematics and
judge mathematicians of the past according to criteria developed within the profession,
and ignore any kind of political fads and pressures. 
\end{abstract}

\title[`Decolonisation' of the curricula]{`Decolonisation' of the curricula\\ and some related issues}

\author{Alexandre Borovik}
\email{alexandre@borovik.net}
\date{05 January 2023}

\maketitle

\section*{Disclaimers}

Views expressed are my own and do not necessarily 
represent the position of my (former) employer, or any other person, 
corporation, organisation, or institution.

I also emphasise that, in this paper, I am  talking only about university level mathematics, not high school or primary school  mathematics.

\section{Introduction}
\label{intro}
I am writing this because my colleagues asked me to comment on a request distributed in one of the British universities:

\begin{quote}
\emph{The majority of the mathematical results in our undergraduate courses are attributed to white European men. This is questionable  because this is not representative of the true history of mathematics. It also sends the message that research in mathematics can be done only by people of the same background. We all have to think how contributions of mathematicians from other backgrounds can be more proportionately represented in our curriculum.}
\end{quote}
This moderate proposal, made, obviously, with the best intentions, belongs to the current trend  of  `\emph{decolonisation of education}'.  It raises a number of questions, which I  will try to outline in this paper.

First of all, in other academic disciplines understanding of `decolonisation' is much wider and includes removal from courses any mentioning of works of authors, and authors themselves (also known as `\emph{cancelling'}), if their behaviour can be judged as  colonialist, racist, misogynist, or detrimental in any way to the aspirations and  rights of any underprivileged or oppressed minority. The demands of `decolonisation' are universal, they are not subject-specific. It is my guess that very soon they will be applied to mathematics in their full strength.

Moreover, mathematics itself is seen in certain quarters (mostly in the US) as a weapon of oppression and intimidation of underprivileged racial and ethnic groups. These ideas are clearly expressed in the \emph{Equitable Math} programme \cite{Pathway2020}. All this debate is about secondary school mathematics teaching. But the tide of changes is already reaching university mathematics as well. In Britain this is becoming the issue of academic freedom and freedom of speech, see, for example, \emph{We are too scared to challenge 'decolonisation' of maths curriculum, top professors warn} \cite{Clarence-Smith2022}, a recent article in a influential newspaper \emph{The Telegraph} -- it is part of debates around the \href{https://bills.parliament.uk/bills/2862}{Higher Education (Freedom of Speech) Bill} which is currently passing through the British Parliament.

I can now formulate the key message of my paper:
\bq
\emph{The international mathematical community should defend professional freedom and insist on our right to formulate our curricula and evaluate the history of mathematics and judge mathematicians of the past according to criteria developed within the profession, and ignore any kind of political fads and pressures.}
\eq

Not long ago I have happily retired  (this is perhaps why my colleagues asked me to comment publicly on this matter). So I decided to run a mental experiment: how would I  `decolonise' my linear algebra course which I taught for more than 20 years if I were requested to do so? I think linear algebra is a representative test case, being one of the core courses of undergraduate mathematics, and a mathematical discipline with a millennia  long historic tradition.

A brief summary of my findings:
\bq
\emph{I think undergraduate linear algebra courses would benefit from inclusion of some historical material -- it is interesting and  can add colour.}
\eq

\begin{quote}
\emph{However, a realistic approach to a historical review of  mathematical curricula  is not likely to meet the  expectations of  political forces which instigate the `decolonise education' drive.}
\end{quote}

Well, for us as mathematicians, there is nothing surprising in facing problems, which, under closer consideration,  have no solution. Remember the problem of trisection of an angle and the problem of quadrature of the circle using only compasses and a  straightedge? It is noteworthy that subsequently all past failures were perfectly explained.

I am not a historian of mathematics, I am a practising mathematician. I am open to constructive criticism and I would love to see alternative `decolonisation' suggestions at least for a linear algebra course.

All `decolonisation' projects would benefit from being made known to the wider mathematical community. Publicity and an open free discussion are the best way to resist outside pressures to engage in virtue signalling at the expense of historical and mathematical truth.

\section{My mental experiment: `decolonisation' of linear algebra}
\label{sec:additional-names}

To decolonise my linear algebra course, I would immediately include in my lectures the following observations:

\begin{itemize}
\item  Apparently one of the first problems (designed for teaching!) which involved solving a system of two linear equations with two unknowns was found on a Sumerian cuneiform tablet from about 2000 BC:
    \begin{quote}
\emph{[It] asks for the areas of two fields whose total area is $1800$ sar,
when the rent for one field is $2$ sil\`{a} of grain per $3$
sar, the rent for the other is $1$ sil\`{a} per $2$ sar, and the
total rent on the first exceeds the other by $500$ sil\`{a}.}  \cite{Grcar2011}
    \end{quote}
\item I think every reader of this paper knows that the words \emph{algorithm}  and \emph{algebra}  are Arabic and stem from the same origin. They are formed from the name of Muhammad ibn Musa al-Khwarizmi (c. 780 -- c.850 AD) and the title of his book \emph{The Compendious Book on Calculation by Completion and Balancing}, c. 813--833 AD; the word \emph{al-jabr} means ``completion''or ``rejoining''. However, he was Persian and born in Khwarazm (hence the name), at that time part of the Arab Caliphate, but today situated in Uzbekistan and called Xorazm in Uzbek.

\item  Seki Takakazu   (Japan, c. 1640--1708 AD) introduced determinants and applied them to solving systems of linear equations---considerably earlier than Gabriel Cramer (1704--1752) did that in Europe in 1750.

\item Seki built his theory on much more ancient Chinese methods of solving systems of simultaneous linear equations by Gaussian elimination developed by Liu Hui (c. 225 -- c. 295 AD) in Chapter 8 of his famous commentaries to \emph{Nine Chapters on the Mathematical Art} (263 AD).

\item The famous work by Liu Hui was a commentary on, and improvement of a much older Chinese mathematical text of  about 200 BC. It represented a mathematical tradition which continued up to 17th century AD; \emph{Nine Chapters} was written as a textbook and served this purpose well in China and neighboring countries. Moreover, it was annotated, with some mistakes (perhaps copying errors) corrected c. 640 AD by a committee led by  Li Chunfeng (602--670 AD).
\end{itemize}
Because of my interest in the history of mathematics,  I knew all that earlier. Preparing this paper, I checked, that all these facts could by easily found using only Wikipedia and \href{https://mathshistory.st-andrews.ac.uk/}{MacTutor}: -- mathematics lecturers are not expected to go deeper than that, they have a lot of other stuff to care about.

I would be happy to add historical comments (clearly marked as non-examinable) to the lecture notes. However, in the current political climate, there is a more serious issue: the addition of above historical facts mentioned before is unlikely to resolve the issue of `decolonisation' of linear algebra, because in addition to `white European males' I simply  introduce males from other major civilisations. Alas, I was unable to point to  `\emph{contributions that historically marginalized people have made to mathematics}' (terminology from \cite{California2021}, see also \cite{Replace2021}).

To address this issue, I have to do some additional (maybe superficial, from the point of view of the professional historian) research. In addition, I have to address a number of closely related issues  which I discuss in the later  sections:
\begin{enumerate}
  \item Why do major civilisations feature so prominently in the history of linear algebra?
  \item Why, in the history of mathematics,  mathematicians were so frequently people of privilege?
  \item Where are women?
\end{enumerate}

\section{Why do major civilisations feature so prominently in the history of linear algebra?}
\label{sec:major-civilisations}

For a very simple reason, which applies to all university level mathematics: results and concepts appearing in it were  sufficiently sophisticated and could emerge only at a higher level of technological development and in the presence of division of labour in production, use of money in economic exchanges, and a developed state / political system with taxation, custom duties etc.. This could be explained even to primary school children: if you have to measure and count stuff, it means that you have a lot of it, that you are rich.

The case of ancient civilisations of Mesopotamia illustrates this point well: all the above developments were occurring there the first time in the history of this  land.

Cuneiform tablets with mathematical problems  come from  scribal schools, some of the oldest, maybe the oldest schools in history. Scribal schools trained scribes -- King's (or state) servants who formed a specific caste in the society's structure: they were literate and numerate. It appears that there was no  further specialisation: in the words of the expert,  Jens H{\o}yrup \cite{Hoyrup2013},
\begin{quote}
\emph{The modern colleagues of the scribe are engineers, accountants and notaries.}
\end{quote}
He adds:
\begin{quote}
\emph{If we really want to find Old Babylonian ``mathematicians'' in an
approximately modern sense, we must look to those who created the
techniques and discovered how to construct problems that were
difficult but could still be solved. For example we may think of the
problem TMS XIX \# 2 \emph{[\dots]}: to find the
sides $l$ and $w$ of a rectangle from its area and from the area of another
rectangle[ \dots] whose length is the diagonal of the first rectangle and whose width is the cube constructed on its
length. This is a problem of the eighth degree. Without systematic
work of theoretical character \emph{[\dots] }it would have been impossible to guess that it was bi-biquadratic (our term of course), and that it can be solved by means of
a cascade of three successive quadratic equations.} (p.108)
\end{quote}

Most problems produced for scribal schools were routine exercises for students. However, high level problems as the one above had also been constructed, and they show a sign of a fierce competitiveness within the caste. H{\o}yrup gives  an explanation: \emph{professional pride} (p.109); perhaps also self-advertising and establishing  authority. Scribes were gaining social status and some professional autonomy, and scribal schools played a special role (H{\o}yrup wrote a fascinating book about this development \cite{Hoyrup1991}; a much more updated and recent study is \cite{Michel-Chelma2020}). It is quite possible  that the ruling elite found it useful to educate their children in scribal schools, which was likely to raise the stakes in the competition between schoolmasters. I found a very suggestive piece of evidence in Andrew George's \emph{Introduction} to his translation of \emph{The Epic of Gilgamesh} \cite{George2003}:
\begin{quote}
\emph{For a short period much of Mesopotamia was again united, this time under the kings of the celebrated Third Dynasty of the southern city of Ur, most famously Shulgi (2094--2074 BC in the conventional terminology). The perfect prince was an intellectual as well as warrior and an athlete, and among his many achievements King Shulgi was particularly proud of his literacy and cultural accomplishments. He had rosy memories of his days at the scribal school, where he boasted that he was the most skilled student in his class.} (p. xvii)
\end{quote}
Why does \emph{The Epic of Gilgamesh}, the oldest piece of written literature known to us, appear in this narrative? Scribes maintained \emph{The Epic}, edited it and copied. It was their political manifesto, a polite message to the rulers that they have to behave reasonably and show some mercy to, and maybe even care of, the people they rule.  \emph{The Epic} was routinely copied by scribal schools' students, as an exercise in writing.

A wry comment by  H{\o}yrup \cite[p. 109]{Hoyrup2013} comes as no surprise:
\begin{quote}
\emph{Even analysis of the cultural function of  ``advanced'' Old Babylonian mathematics may thus teach us something about our own epoch.}
\end{quote}

\section{Why, in the history of mathematics,  were mathematicians so frequently people of privilege?}

For centuries, mathematics sufficiently sophisticated to be studied now at university level, was done by people of privilege and for people of privilege. Why? Because only they had sufficient time for doing mathematics. The most important resource for mathematical research is \emph{uninterrupted} time for thinking.  I have an experimental confirmation of that basic  principle: 18 months of lockdown and retirement became perhaps the most productive time of my life. See my paper \cite{Borovik2020} for an informal discussion of my experience.

Unfortunately,  uninterrupted time for thinking  is increasingly scarce in British universities.

Let us inspect the list of mathematicians mentioned in Section \ref{sec:additional-names}.

I have already talked  about Old Babylonian scribes in the previous   Section
 -- they were obviously a privileged caste.

Muhammad ibn Musa al-Khwarizmi served at the court of Caliph Al-Ma'mun (reigned 786--833 AD) in Baghdad. The famous book was written at the Caliph's suggestion. Al-Khwarizmi  was the head of the House of Wisdom, also known as the Grand Library of Baghdad, arguably the most important intellectual centre of the epoch.

This is what is known about Seki Takakazu:
\begin{quote}
\emph{Seki Takakazu} [\dots ] \emph{spent most of his career serving the Tokugawa family in the Koshu Domain as a bureaucrat whose specialty was accounting, and his social position became that of a shogunate retainer because his feudal lord became the successor (later the sixth shogun Tokugawa Ienobu) to the fifth shogun (hereditary military dictator of Japan).} \cite{Edo}
\end{quote}

Liu Hui: very little is known about him. He lived in the northern Wei kingdom during the 3rd century AD, and was descended from aristocracy. A leading specialist in ancient Chinese mathematics Karine Chemla emphasises that
\begin{quote}
\emph{ Liu’s commentary on \emph{The Nine Chapters} proved the correctness of its algorithms. These proofs are the earliest-known Chinese proofs in the contemporary sense.} \cite{Chemla-Liu-Hui}
\end{quote}
She further says:
\begin{quote}
\emph{A certain philosophical perspective permeates the mathematical work of Liu. He quotes a great variety of ancient philosophical texts, such as the Confucian canons, prominently the \emph{Yijing} (I Ching; Book of Changes); Daoist key texts, such as the \emph{Zhuangzi;} and Mohist texts. Moreover, his commentary regularly echoes contemporary philosophical developments. It can be argued that he considered an algorithm to be that which, in mathematics, embodies the transformations that are at play everywhere in the cosmos -- thus his philosophical reflections on mathematics related to the concept of  ``change'' as a main topic of inquiry in China.} \cite{Chemla-Liu-Hui}
\end{quote}
It is clear that Liu Hui was a deep and powerful thinker by standards of any epoch of human history, and he was cognisant of all the contemporary Chinese science and culture. Obviously, he was not underprivileged.

 Li Chunfeng's life illustrates the role of the state in support of mathematics as well as an early example of academic bureaucracy. \emph{MacTutor }\cite{OConnor-Robertson-on-LiChunfeng2003} says about him:
 \begin{quote}
\emph{ Li rose to become a high-ranking court astronomer and historian, being first appointed to the Imperial Astronomical Bureau in $627$ AD. The reason for his appointment was the need for calendar reform. In ancient China there was a belief that a ruler received his right to rule from heaven. Changing the calendar was seen as one of the duties of the office, establishing the emperor's heavenly link on earth. After a change of ruler, and even more significantly after a change of dynasty, the new Chinese emperor would seek a new official calendar thus establishing a new rule with new celestial influences.}
  \end{quote}
  Later in his life, he became director of the Bureau. He was more prominent as an astronomer and astrologer than a mathematician, but, as a high ranking bureaucrat, was appointed as editor-in-chief for a collection of mathematical treatises called the \emph{The Ten Classics}.

\section{Where are women?}

Sadly, I  could not  find theorems or definitions in my lecture course which could be attributed to women, and checking my collection of undergraduate linear algebra textbooks did not help either.  Critical developments in linear algebra took place from the 18th to first part of the 20th century when women were sadly underrepresented in  mathematics -- indeed, most universities did not accept female students or
professors.  On the specific theme of Gaussian elimination there is an informative paper \cite{Grcar2011} which names only two women, Sister Mary Thomas\`{a} Kempis Kloyda (1896--1977) and  Gertrude Louise Macomber \cite{Macomber1923}  -- but as historians of mathematics and not contributors to its development.

\begin{figure}[h]
\begin{center}
 \includegraphics[width=0.5\textwidth]{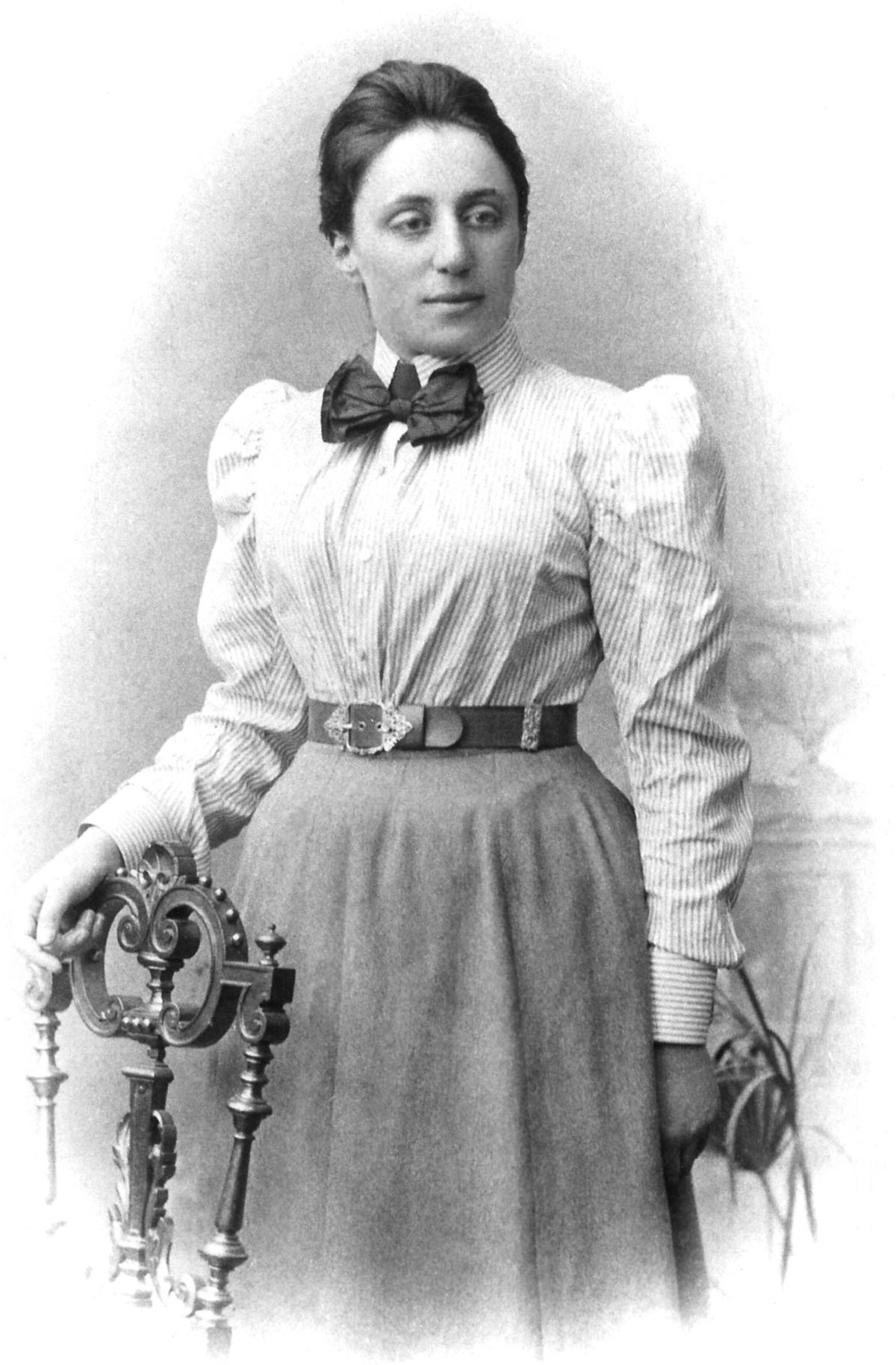}
   \end{center}
\caption{
Emmy Noether (1882--1935).
}
\label{fig:Noether}
\end{figure}

Linear algebra as a taught discipline was born only at the beginning of the 20th century, a period of radical revolution in algebra and the birth of abstract algebra. Again, all the principal contributors were white European males with perhaps a single exception:  Emmy Noether (1882--1935), Figure \ref{fig:Noether}. Noether suffered a direct and grotesque discrimination, both sexist and anti-Semitic, in her academic career, but she was one of the founders of abstract algebra.  Her ideas in algebra definitely influenced the way we (that is, professional mathematicians) think about linear algebra; unfortunately this higher level of thinking cannot be realistically represented in undergraduate linear algebra courses of most British universities. There is a good example of this high level approach: Chapter XIV \emph{Representation of One Endomorphism} of Lang's \emph{Algebra} \cite{Lang} -- but this is a \emph{graduate}, not undergraduate level book.

Linear algebra  was brought to its final shape in the paper of  1948 by Alan Turing \cite{Turing1948}. He emphasised the role of elementary matrices,  introduced  the $LU$-factorisation of matrices and suggested solving systems of simultaneous linear equations $A\vec{x} = \vec{b}$  for non-degenerate square matrices $A$  not by Gaussian elimination, but by inversion of $A$:  $$\vec{x} = A^{-1}\vec{b}$$ -- all that is a standard stuff in modern textbooks. Turing showed in his paper how this approach allows a better control of rounding errors, but his principal  motivation for his paper was visionary:

\begin{quote}
\emph{It seems probable that with the advent of electronic computers it will become standard practice to find the inverse}.
\end{quote}

\begin{figure}[h]
\begin{center}
 \includegraphics[width=0.6\textwidth]{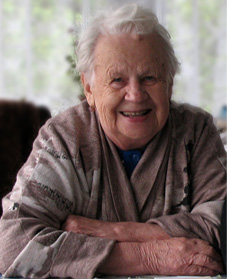}
   \end{center}
\caption{
Vera Nikolaevna Kublanovskaya (1920--2012).
}
\label{fig:Kublanovskaya}
\end{figure}

This indeed has happened, and with the subsequent advance of computers numerical linear algebra was born. In the new socio-economic and political environments after World War II, women became more prominent in numerical  linear algebra than previously  in linear algebra. Since this is out the scope of my course of pure linear algebra, I'll give here just one name; Vera Kublanovskaya (1920--2012), Figure \ref{fig:Kublanovskaya}, who in 1961--62 developed (in parallel with, and independent of, John Francis) what today is called the $QR$-algorithm for finding eigenvalues of a matrix \cite{Golub-Uhlig2009}. I mention her name simply because I remember it from my own days as a student --  I happened to take a course in numerical linear algebra at the time when this  was still  a newborn discipline, and, in particular, computer algorithms for finding eigenvalues of a matrix were a hot topic.

More names can be added, but  I leave this task to my colleagues teaching numerical linear algebra; if they contact me, I'll be happy to help them.

Many mathematical courses carry the birthmarks of development during periods of history with socio-economic and political conditions quite different from those of our times, under conditions of  deeply rooted inequality and grotesque social injustice. But we cannot change history. This has to be accepted as an uncomfortable, but unavoidable fact, and our unfortunate colleagues who teach these courses should not be blamed for the courses' historic origins.

\section{Conclusions}

I do not wish to initiate a debate with the promoters of `decolonisation' of the curricula: I doubt that it would be constructive. In my opinion, the best response to `decolonisation'
pressures is to tell the students the real (pre)history of a particular mathematical discipline -- I call this activity  ``historical enrichment''.  In a few cases, this will bring revelations;
in calculus, for example, the fabulous Kerala mathematics of the 15th--16th centuries. But it is up to lecturers in particular courses to decide  which historical examples they wish to
bring into their lectures and up to what level of detail. Historical enrichment can be easier and less time consuming if teaching materials are widely shared and historic examples are
made known to the wider mathematical community. Publicity and an open discussion are the best way to resist outside pressures to engage in virtue signalling at the expense of
historical and mathematical truth.

Our colleagues are lucky if they  teach courses which are rooted in the 19th or even 20th century mathematics. Indeed, how can one `decolonise'  Complex Analysis?  Or Lie
Algebra? Or Functional Analysis? Or Category Theory? However, even there it may be necessary to prepare polite but firm responses to `decolonisation' requests. Let us return to the
quote which was the starting point of this paper.

\begin{quote}
\emph{The majority of the mathematical results in our undergraduate courses are attributed to white European men. This is questionable  because this is not representative of the
true history of mathematics. It also sends the message that research in mathematics can be done only by people of the same background. We all have to think how contributions of
mathematicians from other backgrounds can be more proportionately represented in our curriculum.}
\end{quote}

But what can we do if this simply reflects, in a particular discipline,  the historic truth? Should we be apologetic? I do not think so. We have to offer our full support to, and
solidarity with our colleagues who express the view that there is nothing in their courses justifying a change.

We also have to understand that   the `decolonisation' issue is much wider than mathematics, and that  in the humanities it is much more acute. As my colleague  Elijah Liflyand
wrote to me about fighting `decolonisation' in mathematics,
\begin{quote}
\emph{This is our red line. This is our last redout.}
\end{quote}

This paper is written for my colleagues, university mathematicians. I hope it might help them to determine their position and calmly  stick to it. But I leave it to them to decide what
is the most appropriate course of action for them. They are mathematicians; in their mathematical work, they are free and independent -- this is the nature of mathematics.

\section*{Acknowledgements}
I thank Gregory Cherlin, Tony Gardiner,  Michael Grinfeld, Mikhail Katz, Alexander Kheyfits, Roman Kossak, Elijah Liflyand, Alexei Muravitsky, Dima Pasechnik, Evgeny Plotkin,  Sergey Shpectorov, and George Wilmers for their criticism, corrections (up to the level of copy-editing!),  and advice. The help  they gave me does not mean that they necessarily endorse all (or any) of the claims that I made in this paper.

I also thank all my colleagues who commented on this paper but preferred to stay anonymous -- and I am deeply saddened to see  that they have good reasons indeed  to avoid any association with this paper.

\end{document}